\documentclass[11pt]{article}
\usepackage{amssymb}

\title{Bijective Proofs for ``Enumerative Properties of Ferrers Graphs''}
\author{Jason Burns}
\begin{document}
\maketitle

 
\section{Introduction}
In their paper \cite{ehrenborg}, Richard Ehrenborg and Stephanie van Willigenburg introduce the concept of a ``Ferrers graph''. 
Given a Ferrers diagram,%
\footnote{Recall that a Ferrers diagram is a set of lattice points in the upper-right quadrant 
\( \{(x,y):x,y\in\mathbb{N}\} \) such that, if $(x,y)$ is in the set, so is $(x',y')$ for all $x'$ at most $x$ and $y'$ at most $y$. (We normally rotate the diagram by $-\frac{\pi}{2}$ and draw boxes for the lattice points.)}
we construct its corresponding Ferrers graph by taking the vertices to be the rows and columns of the diagram, and the edges to be the squares of the diagram: row $a$ and column $b$ are connected by an edge if and only if there is a square in position $(a,b)$ of the diagram.

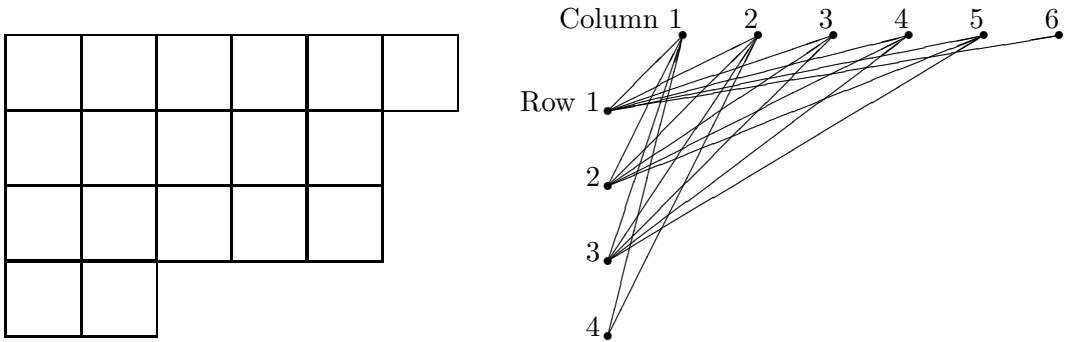
\begin{figure}[htbp] 
\centering
{\setlength{\unitlength}{1.0cm}
\begin{picture}(15,4.5)
	\put(8,3){\line(1,1){1}} 
	\put(8,3){\line(2,1){2}} 
	\put(8,3){\line(3,1){3}}
	\put(8,3){\line(4,1){4}}
	\put(8,3){\line(5,1){5}}
	\put(8,3){\line(6,1){6}}
	\put(8,2){\line(1,2){1}}
	\put(8,2){\line(1,1){2}} 
	\put(8,2){\line(3,2){3}}
	\put(8,2){\line(2,1){4}}
	\put(8,2){\line(5,2){5}}
	\put(8,1){\line(1,3){1}}
	\put(8,1){\line(2,3){2}}
	\put(8,1){\line(1,1){3}}
	\put(8,1){\line(4,3){4}}
	\put(8,1){\line(5,3){5}}
	\put(8,0){\line(1,4){1}}
	\put(8,0){\line(1,2){2}}

	\put(8,3){\circle*{0.1}} 
	\put(8,2){\circle*{0.1}}
	\put(8,1){\circle*{0.1}}
	\put(8,0){\circle*{0.1}}
	\put(9,4){\circle*{0.1}}
	\put(10,4){\circle*{0.1}}
	\put(11,4){\circle*{0.1}}
	\put(12,4){\circle*{0.1}}
	\put(13,4){\circle*{0.1}}
	\put(14,4){\circle*{0.1}}

	\put(7.9,3){\makebox(0,0)[br]{Row 1}} 
	\put(7.9,2){\makebox(0,0)[br]{2}} 
	\put(7.9,1){\makebox(0,0)[br]{3}}
	\put(7.9,0){\makebox(0,0)[br]{4}}
	\put(9,4.1){\makebox(0,0)[br]{Column 1}}
	\put(10,4.1){\makebox(0,0)[br]{2}}
	\put(11,4.1){\makebox(0,0)[br]{3}}
	\put(12,4.1){\makebox(0,0)[br]{4}}
	\put(13,4.1){\makebox(0,0)[br]{5}}
	\put(14,4.1){\makebox(0,0)[br]{6}}
	
	\put(0,3){\framebox(1,1){ }}
	\put(1,3){\framebox(1,1){ }}
	\put(2,3){\framebox(1,1){ }}
	\put(3,3){\framebox(1,1){ }}
	\put(4,3){\framebox(1,1){ }}
	\put(5,3){\framebox(1,1){ }}
	\put(0,2){\framebox(1,1){ }}
	\put(1,2){\framebox(1,1){ }}
	\put(2,2){\framebox(1,1){ }}
	\put(3,2){\framebox(1,1){ }}
	\put(4,2){\framebox(1,1){ }}
	\put(0,1){\framebox(1,1){ }}
	\put(1,1){\framebox(1,1){ }}
	\put(2,1){\framebox(1,1){ }}
	\put(3,1){\framebox(1,1){ }}
	\put(4,1){\framebox(1,1){ }}
	\put(0,0){\framebox(1,1){ }}
	\put(1,0){\framebox(1,1){ }}
\end{picture}
}
\caption{A Ferrers diagram and its graph.}
\end{figure}

This bipartite graph is clearly equivalent to the original diagram, but in this new form we can raise new questions.  Ehrenborg and Van Willigenburg examine several of these, including the number of spanning trees and Hamiltonian paths of the Ferrers graph; alas, the proofs they give for these are not bijective.  In this paper, I give bijective proofs for both of these, which extend to weighted versions.

Ehrenborg and Jeffrey B. Remmel have independently found a bijection for the number of spanning trees; I have not yet seen it.

\section{Hamiltonian paths}

Let's suppose our Ferrers diagram has the same number $n$ of rows and columns.%
\footnote{This is not much of a restriction. Since any Hamiltonian path in the Ferrers graph must alternate between rows and columns, there can be at most one more (or fewer) row than column, if there are to be any Hamiltonian paths at all.}
Ehrenborg and Van Willigenburg proved that the number of Hamiltonian paths for the corresponding Ferrers graph is the square of the number of $n$--rook placements on the Ferrers diagram. (By ``number of $n$--rook placements'', we mean the number of ways to choose $n$ squares from the diagram so that no two are in the same row or column.) To prove this bijectively, we seek a procedure converting any Hamiltonian path to a pair of $n$--rook placements, and an inverse procedure converting any pair of
$n$-rook placements back to a Hamiltonian path.

\begin{figure}[htbp]
\centering
{\setlength{\unitlength}{1.0cm}
\begin{picture}(7,5.5)
	{\thinlines 
	\put(1,4){\line(1,1){1}}
	\put(1,4){\line(2,1){2}}
	\put(1,4){\line(3,1){3}}
	\put(1,3){\line(1,1){2}}
	\put(1,3){\line(2,1){4}}
	\put(1,2){\line(2,3){2}}
	\put(1,2){\line(1,1){3}}
	\put(1,1){\line(1,4){1}}
	\put(1,0){\line(1,5){1}}
	}
	{\thicklines 
	\put(1,4){\line(4,1){4}}
	\put(1,4){\line(5,1){5}}
	\put(1,3){\line(1,2){1}}
	\put(1,3){\line(3,2){3}}
	\put(1,2){\line(1,3){1}}
	\put(1,2){\line(4,3){4}}
	\put(1,1){\line(1,2){2}}
	\put(1,1){\line(3,4){3}}
	\put(1,0){\line(2,5){2}}
	}

	\put(1,4){\circle*{0.1}}
	\put(1,3){\circle*{0.1}}
	\put(1,2){\circle*{0.1}}
	\put(1,1){\circle*{0.1}}
	\put(1,0){\circle*{0.1}}
	\put(2,5){\circle*{0.1}}
	\put(3,5){\circle*{0.1}}
	\put(4,5){\circle*{0.1}}
	\put(5,5){\circle*{0.1}}
	\put(6,5){\circle*{0.1}}

	\put(0.9,4){\makebox(0,0)[br]{Row 1}} 
	\put(0.9,3){\makebox(0,0)[br]{2}} 
	\put(0.9,2){\makebox(0,0)[br]{3}}
	\put(0.9,1){\makebox(0,0)[br]{4}}
	\put(0.9,0){\makebox(0,0)[br]{5}}
	\put(2,5.1){\makebox(0,0)[br]{Column 1}}
	\put(3,5.1){\makebox(0,0)[br]{2}}
	\put(4,5.1){\makebox(0,0)[br]{3}}
	\put(5,5.1){\makebox(0,0)[br]{4}}
	\put(6,5.1){\makebox(0,0)[br]{5}}
\end{picture}
}
\caption{A Hamiltonian path in a Ferrers graph.}
\end{figure}
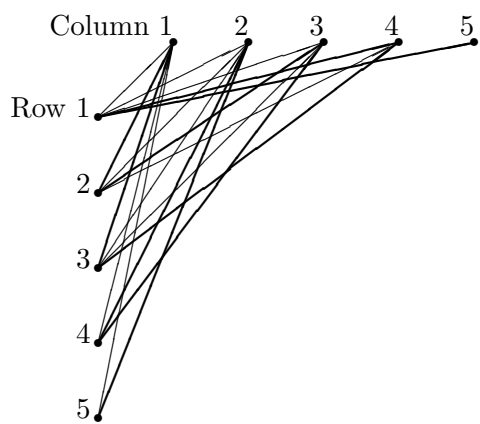


\emph{Hamiltonian paths to pairs of $n$--rook placements}

One end of the Hamiltonian path is a row-vertex, and the other is a column-vertex. Suppose the vertices have labels $a_1, b_1, a_2, b_2, \dots, a_n, b_n$, starting at the row end, so the ``$a$'' vertices are rows and the ``$b$'' vertices are columns.

\begin{enumerate}
\item Start with an empty Ferrers diagram which we'll mark up with two distinct $n$--rook placements, labeled with A's and B's. Also set \(b_0 := 1\).
\item Repeat for each pair of vertices $a_i, b_i$:
	\begin{enumerate}
	\item Mark an A in the square \( (a_i,b_{i-1}) \), unless there is already an A in  column $b_{i-1}$.
If there is an A in that column, mark the A in the next available column without an A.
	\item Mark a B in the square \( (a_i, b_i) \).
	\item Move on to the next pair of vertices $a_{i+1}, b_{i+1}$.
	\end{enumerate}
\end{enumerate}

\begin{figure}[htbp]
\centering
{\setlength{\unitlength}{1.0cm}
\begin{picture}(5,5.5)
	\put(0,4){\framebox(1,1){ }} 
	\put(1,4){\framebox(1,1){ }}
	\put(2,4){\framebox(1,1){ }}
	\put(3,4){\framebox(1,1){ }}
	\put(4,4){\framebox(1,1){AB}} 
	\put(0,3){\framebox(1,1){B}}
	\put(1,3){\framebox(1,1){ }}
	\put(2,3){\framebox(1,1){A}}
	\put(3,3){\framebox(1,1){ }}
	\put(0,2){\framebox(1,1){ }}
	\put(1,2){\framebox(1,1){ }}
	\put(2,2){\framebox(1,1){ }}
	\put(3,2){\framebox(1,1){AB}}
	\put(0,1){\framebox(1,1){ }}
	\put(1,1){\framebox(1,1){A}}
	\put(2,1){\framebox(1,1){B}}
	\put(0,0){\framebox(1,1){A}}
	\put(1,0){\framebox(1,1){B}}
\end{picture}
}
\caption{The $n$--rook configuration corresponding to the Hamiltonian path above.}
\end{figure}

We place one B in each row $a_i$, and one B in each column $b_i$, so the B's constitute an $n$--rook placement. We place one A in each row $a_i$, and we never place more than one A in any column, so the A's also constitute an $n$--rook placement. The only thing we need to check is that there always is a ``next available column'' in row $a_i$, in which to place our A.

This seems easy: after all, we're about to mark a B in row $a_i$, column $b_i$, which guarantees that column $b_i$ is possible. Since we mark as many A's as B's at every stage, and every column we've marked so far with a B (that is, $b_1$ through $b_{i-1}$) is guaranteed to have an A, we also know that $b_i$ doesn't have an A and is thus available.  So there is at least one available column.

The problem is that we need a \emph{next} available column, that is, one that follows $b_{i-1}$.
The solution is to observe that the ``next available column'' is actually the \emph{first} available column:
every column above it has an A already. Suppose, at some pair of vertices $a_i, b_i$, that the column $b_{i-1}$ we would have used already has an A. We've never been to $b_i$ before, so there are only two ways $b_{i-1}$ could have acquired that A. Either it was put there in the first iteration of the loop,
and $b_{i-1}=b_0$ is the first column (in which case the next column without an A is obviously the first one without an A), or else it was put there as the ``next available column'' (so $b_{i-1}$ was at one time the first available column, which means that all the columns before $b_{i-1}$ contain an A).
Hence the available column $b_i$ must follow $b_{i-1}$.


\emph{Pairs of $n$--rook placements to Hamiltonian paths}

\begin{enumerate}
\item Start with a Ferrers diagram marked with two different $n$--rook placements, labeled with A's and B's. Also set \(b_0 := 1\).
\item Repeat until we've identified each pair of vertices $a_i, b_i$:
	\begin{enumerate}
	\item Remove the A in column $b_{i-1}$, or else the first column with an A, and set $a_i$ equal to the row it was in.
	\item Remove the B in row $a_i$, and set $b_i$ equal to the column it was in.
	\item Move on to identify the next pair of vertices $a_{i+1}, b_{i+1}$.
	\end{enumerate}
\end{enumerate}

This is clearly the inverse of the procedure given earlier. Note that ``first column with an A'' is equivalent to ``next column with an A'', by our earlier argument.

One last comment: if we switch the roles of rows and columns in the above procedures, or if we switch the roles of A's and B's, we get another Hamiltonian path. It would be interesting to know if there is any deeper meaning to these two involutions.

\section{Spanning trees}

The number of spanning trees of a Ferrers graph is the product of the lengths of all but the first column, times the product of the lengths of all but the first row.  We can prove this with a similar bijection.

We represent this product combinatorially in a similar fashion to the rook placements above: in each row but the first we place an R, and in each row but the first we place a C. (Note that now there may be more than one R in some column, or more than one C in some row; and we needn't have equal number of rows and columns.) Also, place an X in the first row and first column, for a reason we'll see later.

\begin{figure}[htbp]
\centering
{\setlength{\unitlength}{1.0cm}
\begin{picture}(5,4.5)
	\put(0,3){\framebox(1,1){X}}
	\put(1,3){\framebox(1,1){ }}
	\put(2,3){\framebox(1,1){ }}
	\put(3,3){\framebox(1,1){C}}
	\put(4,3){\framebox(1,1){ }}
	\put(5,3){\framebox(1,1){C}}
	\put(0,2){\framebox(1,1){ }}
	\put(1,2){\framebox(1,1){ }}
	\put(2,2){\framebox(1,1){RC}}
	\put(3,2){\framebox(1,1){ }}
	\put(4,2){\framebox(1,1){ }}
	\put(0,1){\framebox(1,1){ }}
	\put(1,1){\framebox(1,1){R}}
	\put(2,1){\framebox(1,1){ }}
	\put(3,1){\framebox(1,1){ }}
	\put(4,1){\framebox(1,1){C}}
	\put(0,0){\framebox(1,1){ }}
	\put(1,0){\framebox(1,1){RC}}
\end{picture}
}
\caption{A configuration corresponding to a spanning tree.}
\end{figure}

Given such a configuration of R's and C's, we will convert it to a spanning tree in two stages.

For the first stage, note that every row except the first has exactly one R, but it need not have any C's; and likewise, each column but the first has one C and as few as zero R's.  If there is a row $a$ with an R in column $b$, but no other markings, we may construct a spanning tree by removing row $a$, constructing the spanning tree on the remaining rows and columns, and then joining row $a$ to the rest of the tree with a single edge to column $b$. Likewise, if there is a column $b$ with a C in column $a$, but no other markings, we may construct a spanning tree by removing column $b$, constructing the spanning tree on the remaining rows and columns, and then joining column $b$ to the rest of the tree
with a single edge to row $a$. It clearly doesn't matter in which order we prune these rows and columns.  Conversely, given a spanning tree, we can remove a leaf node (say, row-vertex $a$), mark the square $(a,b)$ corresponding to the edge connecting it with the rest of the tree (we would mark with an R, since this is a row-vertex), and proceed to mark the remaining rows and columns of the spanning tree. Thus far, then, we have a bijection. Note that we cannot remove the first row or column in this stage.

For the second stage, we are thus left with only the irreducible objects of each type: the configurations of R's and C's for which every row but the first has at least one C and every column but the first has at least one R, against the paths from the first row-vertex to the first column-vertex. Remove the first row and column: then what we have left are simply pairs of $n$--rook placements (in the former case), and arbitrary Hamiltonian paths (in the latter), and we already have a bijection between these!
If each row has a C then there must be at least as many C's as R's, and if each column has an R then there are as many R's as C's, and so there are exactly as many R's as C's: one per column, and also one per row. (So the R's are an $n$--rook placement, and so are the C's.)
On the other hand, a Hamiltonian path starting at the first column and ending at the first row is simply an arbitrary Hamiltonian path on the second through $n$th rows and columns; just strip off the first-row and first-column endpoints, or put them back on again.

\begin{figure}[htbp]
\centering
{\setlength{\unitlength}{1.0cm}
\begin{picture}(8,4.5)
	{\thinlines 
	\put(1,3){\line(1,1){1}} 
	\put(1,3){\line(2,1){2}} 
	\put(1,3){\line(5,1){5}}
	\put(1,2){\line(1,2){1}}
	\put(1,2){\line(2,1){4}}
	\put(1,2){\line(5,2){5}}
	\put(1,1){\line(1,3){1}}
	\put(1,1){\line(1,1){3}}
	\put(1,1){\line(4,3){4}}
	}
	
	{\thicklines 
	\put(1,3){\line(3,1){3}}
	\put(1,3){\line(4,1){4}}
	\put(1,3){\line(6,1){6}}
	\put(1,2){\line(1,1){2}} 
	\put(1,2){\line(3,2){3}}
	\put(1,1){\line(2,3){2}}
	\put(1,1){\line(5,3){5}}
	\put(1,0){\line(1,4){1}}
	\put(1,0){\line(1,2){2}}
	}

	\put(1,3){\circle*{0.1}} 
	\put(1,2){\circle*{0.1}}
	\put(1,1){\circle*{0.1}}
	\put(1,0){\circle*{0.1}}
	\put(2,4){\circle*{0.1}}
	\put(3,4){\circle*{0.1}}
	\put(4,4){\circle*{0.1}}
	\put(5,4){\circle*{0.1}}
	\put(6,4){\circle*{0.1}}
	\put(7,4){\circle*{0.1}}

	\put(0.9,3){\makebox(0,0)[br]{Row 1}} 
	\put(0.9,2){\makebox(0,0)[br]{2}} 
	\put(0.9,1){\makebox(0,0)[br]{3}}
	\put(0.9,0){\makebox(0,0)[br]{4}}
	\put(2,4.1){\makebox(0,0)[br]{Column 1}}
	\put(3,4.1){\makebox(0,0)[br]{2}}
	\put(4,4.1){\makebox(0,0)[br]{3}}
	\put(5,4.1){\makebox(0,0)[br]{4}}
	\put(6,4.1){\makebox(0,0)[br]{5}}
	\put(7,4.1){\makebox(0,0)[br]{6}}
\end{picture}
}
\caption{The spanning tree.}
\end{figure}
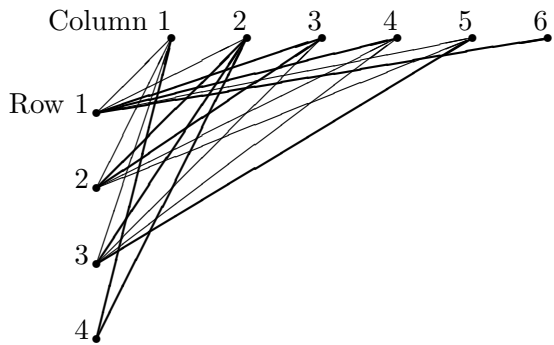

\section{Weighted spanning trees}

Our spanning-tree bijection also gives us a weighted formula for the total number of spanning trees.
Specifically, if we assign weight $x_a y_b$ to the edge from row $a$ to column $b$, then the weight of any given spanning tree is equal to the weight of the corresponding configuration of R's, C's, and one X in the upper-left corner.
This in turn lets us sum the weights of all spanning trees by summing the weights of all configurations, giving the product formula
\[
(x_1 x_2 \cdots x_m)\cdot(y_1 y_2 \cdots y_n) \cdot
\prod_{2\leq a\leq m}(y_1 + \cdots + y_{\lambda_a^{ }})
\cdot \prod_{2\leq b \leq n}(x_1 + \cdots + x_{\lambda_b^\prime})
\]
where $\lambda_a$ and $\lambda_b^\prime$ are the lengths of row $a$ and column $b$ respectively.%
\footnote{This is not hard to see. Since each symbol in the configuration can be placed independently of the others, we simply calculate the weights each symbol can give, and multiply them together. The X must be placed in the upper left corner, contributing $x_1 y_1$; the R in row $a$ can be placed in column \(1, 2, \dots, \lambda_a\), for a factor of \(x_a\cdot (y_1+\cdots+y_{\lambda_a})\); and the C in row $b$ can be placed in rows \(1,2,\dots,\lambda_b^\prime\), for a factor of \(y_b \cdot (x_1+\cdots+x_{\lambda_b^\prime})\). }

That our spanning-tree bijection is weight-preserving is evident in the first stage, because we remove an edge of weight $x_a y_b$ from both the configuration and the spanning tree; in the second stage it is only slightly less obvious.
Each Hamiltonian path in the second stage starts in the first column and ends in the first row, but passes through all other vertices; hence it contributes \(x_1 y_1 (x_2\cdots x_m)^2 (y_2\cdots y_n)^2\) to the total. Likewise, there is a single X in the first row and the first column, but an R and a C in each other row and column, so these marks also contribute \(x_1 y_1 (x_2\cdots x_m)^2 (y_2\cdots y_n)^2\). 

\section{Future research}

One natural extension%
\footnote{Suggested by Ehrenborg and Van Willenburg in their paper.}
of the definition of a Ferrers graphs is to skew Ferrers diagrams.
In fact, the definition works with no further changes; unfortunately, the results above fail, and no reasonable fix is apparent.
Nothing interesting is yet known in this case.

Let me suggest another extension which may prove more fruitful.  One of the key invariants of a graph is its Tutte polynomial, which can be specialized to give such quantities as the number of spanning trees and the chromatic polynomial.
(Thus, a combinatorial determination of the Tutte polynomial would encompass the result above.%
\footnote{Ehrenborg and Van Willenburg also suggested investigating the Tutte polynomial.}
)
The Tutte polynomial can be computed recursively by deleting and contracting edges, but if we try to calculate the Tutte polynomial of a Ferrers graph recursively by removing rows, these contractions yield multiple edges.  Specifically, we get what we might call ``plane-partition graphs'': the number of edges from row $a$ to column $b$ (corresponding to square $(a,b)$) does not increase as we increase $a$ or $b$ --- or, stated another way, if we were to write in square $(a,b)$ the number of edges corresponding to that square, we would get a plane partition.

\end{document}